\documentclass{article}
\usepackage[utf8]{inputenc}
\usepackage{hyperref}
\usepackage{tikz,lipsum,lmodern}
\usepackage[normalem]{ulem}
\makeatother
\usepackage[most]{tcolorbox}
\usepackage{amsmath}
\usepackage{amssymb}
\usepackage{amsthm}
\usepackage{amsfonts}
\usepackage{fancyhdr}
\usepackage{indentfirst}
\usepackage[pages=some]{background}
\usepackage{tikzpagenodes}
\usepackage{lipsum}
\usepackage[left=1cm, right=2cm]{geometry}
\backgroundsetup{
  scale=1,angle=0,color=black,
  contents={%
  \tikz[remember picture,overlay]{
    ([yshift=120pt,xshift=5pt]current page.south west) 
    rectangle 
    ([yshift=5pt,xshift=-5pt]current page.south east);
    \node[anchor=north west,inner sep=0pt,text width=.3\paperwidth,align=left] at ([yshift=-20pt,xshift=10pt]current page.south west|-current page text area.south west) {Konstantinos Konstantinidis \\ Thessaloniki, Greece\\ \href{mailto: kkonsa@ece.auth.gr}{ kkonsa@ece.auth.gr}};
  }
  }
}
\usepackage[left=1cm, right=2cm]{geometry}
\backgroundsetup{
  scale=1,angle=0,color=black,
  contents={%
  \tikz[remember picture,overlay]{
    ([yshift=120pt,xshift=5pt]current page.south west) 
    rectangle 
    ([yshift=5pt,xshift=-5pt]current page.south east);
    \node[anchor=north west,inner sep=0pt,text width=.3\paperwidth,align=left] at ([yshift=-20pt,xshift=10pt]current page.south west|-current page text area.south west) {Konstantinos Konstantinidis \\ Thessaloniki, Greece\\ \href{mailto: kkonsa@ece.auth.gr}{ kkonsa@ece.auth.gr}};
  }
  }
}

\title{Real Analysis in Functional Equations}
\author{Konstantinos Konstantinidis  \\  \\ \href{mailto: kkonsa@ece.auth.gr}{kkonsa@ece.auth.gr}}
\date{11 November 2023}

\begin{document}

\maketitle
\section*{\color{blue}Acknowledgements} I express my sincere gratitude to my teacher, \textbf{Athanasios Kontogeorgis}, for his meticulous review of this paper and to \textbf{Minos Margaritis} for his invaluable suggestions.
\section*{\color{blue}Contents}
\hyperref[test]{\color{blue} 1. Introduction \hfill 2}\\
\indent\hyperref[test1]{\color{blue} 2. Squeeze theorem and Limits of Sequences\hfill 3}\\
\indent\hyperref[test2]{\color{blue} 3. Supremum and Infimum \hfill 6}\\
\indent\hyperref[test3]{\color{blue} 4. Limits Chasing \hfill 14}\\
\indent\hyperref[test4]{\color{blue} 5. Differentiability \hfill 19}\\
\indent\hyperref[test5]{\color{blue} 6. Problems for Practice \hfill 23}\\
\indent\hyperref[test6]{\color{blue} 7. References \hfill 24}\\
\pagebreak
\\
\\
\begin{center}
\textbf{Abstract}\\
In this article, we will showcase some analytical concepts that can be used to tackle Functional Equations (FE) in the positive real numbers domain. Such concepts and related techniques have occasionally appeared in recent High School Math Olympiads, and they are often accompanied neatly by other known techniques. In each section, we develop a theoretical background; next, we briefly mention methods that employ the theory and conclude the article by providing unsolved problems that the reader can try independently.\\
\end{center}

\pagebreak
\phantomsection\label{test}{\color{blue}\section{Introduction}}
Essentially, in this article we intrinsically employ two powerful methods: \textbf{Bounding} and \textbf{Optimizing}. In rough terms, in every example the goal is to find two "almost equal" functions $f_{n}$ and $g_{n}$ wrapping the function we are looking for, i.e, \\
\[f_{n}(x)\leq f(x)\leq g_{n}(x)\quad\forall{n\in\mathbb{N}}\]
Then, by letting $n\to\infty$ we can derive the closed form of $f$.\\
\indent In order to "wrap" the function $f$, we try to "optimize" the bounds found on the $n-1$ step by bounding a side (either the left or the right) of the given functional by imposing the inequality of the $n-1$ step to each term. We simplify and hope that $f_{n}$ would be a better approximation for $f$ than $f_{n-1}$.\\
\indent Usually, the simply stated "simplification" can be cumbersome and would be too technical to perform. This is when the supremum and infimum come in handy. More specifically, rather than following a "bottom-up" approach, i.e., starting from $f_{0},g_{0}$ and generating $f_{n}, g_{n}$, we follow a "top-down" approach, i.e., we assume that we have already found the functions $\sup f\; (x)=f_{n}(x)$ \hyperref[test00]{\color{blue}$^{1}$} and $\inf f\;(x)=g_{n}(x)$ that best approximate $f$ and proceed by contradiction. In other words, any other approximation for $f$ should be worse than the one given by $\sup f\; (x)$, $\inf f\; (x)$.\\ 
\indent Furthermore, instead of "globally" applying this method, i.e., searching for $f_{n}, g_{n}$ such that all $x$ obey the inequality, we can apply it "locally". In other words, we can search for $f_{n}, g_{n}$ such that:\\
\[f_{n}(x_{0})<f(x_{0})<g_{n}(x_{0})\]
for a particular $x_{0}$ and derive the limit $l=\lim_{x\to x_{0}}f(x)$ and use it to find $f$. In particular, when $x_{0}$ is a boundary point (e.g., $x_{0}=0$ for positive real-valued functions), the monotonicity of $f$ can connect the limit $l$ with the $\sup$ or the $\inf$ of $f$.\\
\indent Usually, the solution we expect is of the form $f(x)=cx$, therefore, it is advisable to apply this method to $\displaystyle \frac{f(x)}{x}$ and constants $f_{n}(x), g_{n}(x)$ to deduce that $\frac{f(x)}{x}$ or to calculate the limit $\lim_{x\to x_{0}}\frac{f(x)}{x}$, i.e, the derivative of $f$ at $x_{0}$. Though at first glance the limit $\lim_{x\to x_{0}}\frac{f(x)}{x}$ might not be as useful as the $\lim_{x\to 
x_{0}}f(x)$ we have included in this article some examples that illustrate methods to harness the existence of the derivative to solve the problem, such as the derivation of a differential equation.\\
\indent The concepts briefly mentioned previously and further explained in the following sections provide a broad toolbox for solving even intricate functional equations problems. I hope, therefore, that this article will not only communicate to the reader this toolbox but also bestow valuable insights for solving such problems.
\\
\\
\\
\\
\\
\\
\\
\\
\\
\\
\\
\\
\\
\\
\\
\begin{scriptsize} \phantomsection\label{test00}{\color{blue}$^{1}$} only for the purposes of the introduction, by abuse of notation we assumed that "$\sup f \;(x)$" ($\neq \sup f(x))$ is a function rather than number. We mean that" $\sup f \;(x)$" is the function such that $\sup f \;(x)>f(x)$ and $\forall n\;\in\mathbb{N}$, $f_{n}>\sup f \;(x)$
\end{scriptsize}
\pagebreak
\phantomsection\label{test1}{\color{blue}\section{Squeeze Theorem and Limits of Sequences}}
\begin{tcolorbox}[colback=red!5!white,colframe=red!75!black]
\textbf{Definition 2.1}: We call $x$ the limit of the sequence $x_{n}$ if the following condition holds:
For each real number $\varepsilon > 0$, there exists a natural number $N$ such that, for every natural number $n\geq N$ we have: \[|x_{n}-x|<\varepsilon\]
\end{tcolorbox}
Some properties of limits of real sequences that we will be using extensively later are the following (provided that the limits on the right exist):\\
\begin{tcolorbox}[colback=blue!5!white,colframe=blue!75!black,title=Properties of limits]
\begin{itemize}
\item The limit of a sequence is unique.
\item $ \displaystyle \lim _{n\to \infty }(a_{n}\pm b_{n})=\lim _{n\to \infty }a_{n}\pm \lim _{n\to \infty }b_{n}$
\item $\displaystyle\lim _{n\to \infty }ca_{n}=c\cdot \lim _{n\to \infty }a_{n}$
\item $\displaystyle\lim _{n\to \infty }(a_{n}\cdot b_{n})=(\lim _{n\to \infty }a_{n})\cdot (\lim _{n\to \infty }b_{n})$
$\displaystyle\lim_{n\to \infty }\left({\frac {a_{n}}{b_{n}}}\right)=\frac {\lim \limits _{n\to \infty }a_{n}}{\lim _{n\to \infty }b_{n}}$
provided $\lim_{n\to \infty }b_{n}\neq 0$
\item $\displaystyle\lim _{n\to \infty }a_{n}^{p}=\left[\lim _{n\to \infty }a_{n}\right]^{p}$
\item If  $a_{n}\leq b_{n}$ for all $n$ greater than some $N$, then ${\displaystyle \lim _{n\to \infty }a_{n}\leq \lim _{n\to \infty }b_{n}.}$
\item \textbf{Squeeze Theorem:} If $a_{n}\leq c_{n}\leq b_{n}$ for all  $n>N$ and $\displaystyle\lim _{n\to \infty }a_{n}=\lim _{n\to \infty }b_{n}=L,$ then $\displaystyle \lim _{n\to \infty }c_{n}=L.$
\item If a sequence is bounded and monotonic, it is convergent.
\item A sequence is convergent if and only if every subsequence is convergent.
\item If every subsequence of a sequence has its own subsequence that converges to the same point, then the original sequence converges to that point.
\item \textbf{Bolzano-Weirstrass Theorem:} Every bounded sequence has at least one convergent subsequence.
\item If $\displaystyle\lim_{n\to\infty} a_{n}\to {l}$, f is a continuous function, and $a_{n+1}=f(a_{n})$, then $f(l)=l$
\end{itemize}
\end{tcolorbox}

\begin{tcolorbox}[enhanced,fit to height=2.8cm,
  colback=green!25!black!10!white,colframe=green!75!black,title=Squeezing Technique,
  drop fuzzy shadow,watermark color=white]
    When we have guessed that the function $f$ we are looking for has the form $f(x)=cx$, we often try to find two sequences $(a_{n}),(b_{n})$ converging to the same limit such that $a_{n}x\geq f(x) \geq b_{n}x$ for all n and an arbitrary but fixed x. By letting $n\to \infty$ yields $f(x)=cx$ where c is the common limit of $(a_{n}),(b_{n})$.\\
\end{tcolorbox}
Now let's see an example:\\
\begin{tcolorbox}[colback=green!5!white,colframe=green!75!black]
Problem 2.1
  \tcblower
Find all functions $f:\mathbb{R^{+}}\to \mathbb{R^{+}}$ such that:
$f(x)+f(f(x))=2x$ for all $x\in\mathbb{R^{+}}$.
\end{tcolorbox}
Though this problem is prone to another (more standard) method by recurrence relations (see Example 11, Continuation 1, [3]), we give another equivalent approach as a first sight to this technique.\\

\textit{Solution:}
 Fix a real number x. We obviously have that $\displaystyle 2x=f(x)+f(f(x))>f(x)$, hence
\begin{align}
f(x)<2x\label{}
\end{align}\\
Plugging $f(x)$ in (1) we obtain: $f(f(x))<2f(x)$.
Hence: $f(x)+f(f(x))=2x< f(x)+2f(x)=3f(x)$ hence $f(x)> \frac{2x}{3}$ (3)
\begin{align}
f(x)> \frac{2x}{3}\label{}
\end{align}
Iterating (2), we deduce that $\displaystyle f(f(x))> \frac{2f(x)}{3}$, and, using (1), we deduce:
$\displaystyle 2x=f(x)+f(f(x)) > f(x)+\frac{2f(x)}{3}=\frac{5f(x)}{3}$ hence \begin{align} f(x)< \frac{6x}{5}\end{align}
We obviously see a pattern, and now we need to formalize it.
We need to find those two sequences $(a_{n}),(b_{n})$ such that $b_{n}x< f(x)< a_{n}x$. To this end, we employ the previous logic with $a_{n}$ instead of fixed numbers (by the previous discussion we can observe that $a_{1}=2$), to specify those sequences.\\
If $f(x)<a_{n}x$ then $f(f(x))< a_{n}f(x)$ hence: \[f(f(x))+f(x)=2x<a_{n}f(x)+f(x)\]
This yields \[ f(x)> \frac{2x}{a_{n}+1}\forall x\] We set $\displaystyle b_{n}=\frac{2}{a_{n}+1}$\hfill (4)\\
Therefore, \[f(f(x))+f(x)=2x> b_{n}f(x)+f(x)\] hence it holds that: $\displaystyle\frac{f(x)}{x}< \frac{2}{b_{n}+1}$. Set $\displaystyle a_{n+1}=\frac{2}{b_{n}+1}$ or (by (4)) $\displaystyle a_{n+1}=\frac{2(a_{n}+1)}{a_{n}+3}$. Note that: $\displaystyle a_{n+1}=2-\frac{4}{a_{n}+3}$\hfill (5).\\
\indent Now, using (5), by induction, one can prove that $a_{n}\geq 1$ for all n. Now, one can easily check that $\displaystyle a_{n+1}=2-\frac{4}{a_{n}+3}\leq a_{n}$ hence $(a_{n})$ is (possibly weakly) decreasing. Since $(a_{n})$ is a bounded decreasing sequence of real numbers, it should converge. \\
\indent Hence, $(a_{n})$ has a limit. To calculate it, we use the last property of limits, i.e., we solve $\displaystyle l=2-\frac{4}{l+3}$ which yields $l=1$ or $l=-2$, with the last case being vacuous, since $f(x)\geq 0$ for all x.
Hence, since $f(x)\leq a_{n}x$ for all $n$, tending $n\to \infty$ we obtain $f(x)\leq x$.\\
\indent We can also work the same way with $b_{n}$, but we may omit this step using the following shortcut:\\
$f(x)\leq x$ hence $f(f(x))\leq f(x)$, so $f(f(x))+f(x)=2x\leq 2f(x)$ so $f(x)\geq x$, hence $f(x)=x$ as desired.\\
\\
\begin{tcolorbox}[colback=green!5!white,colframe=green!75!black]
Problem 2.2 - "A problem that solves itself"
  \tcblower
Find all functions $f:[0,\infty)\to[0,\infty)$ satisfying the functional equation
$f(f(x)-x)=2x\forall x\in[0,\infty)$
\end{tcolorbox}
\textit{Solution:}
In order for LHS to be defined, we get $f(x)\geq x$ $\forall x$. Hence, $f(f(x)-x)\geq f(x)-x$ $\forall x$ $\iff$ $f(x)\le 3x$, thus $x\le f(x)\le 3x$

To construct the sequences mentioned in the previous example, consider $a_nx\le f(x)\le b_nx$, to get $a_n(f(x)-x)\le 2x\le b_n(f(x)-x)$ and so: \[\frac{b_n+2}{b_n}x\le f(x)\le \frac{a_n+2}{a_n}x\]

Hence, the sequences we are looking for are:
\[
a_1=1, \quad
b_1=3,\quad
\displaystyle a_{n+1}=\frac{b_n+2}{b_n}\;\text{(1),} \quad
\displaystyle b_{n+1}=\frac{a_n+2}{a_n}\; \text{(2)}
\]
We see that $a_{1}=1<a_{2}=\frac{5}{3}$. Assume that for all $i\leq k$ the relation $a_{i}<a_{i+1}$ holds.\\
\indent By (2) we have that $b_{i}<b_{i-1}$ for all $i\leq k$, hence, by (1), $a_{k+1}>a_{k}$.
Thus, by induction, $a_n$ is strictly increasing, and $b_n$ is a strictly decreasing sequence of positive real numbers.\\
\indent We can manually check that the limit of $(a_{n})$ is 2, and the limit of $(b_{n})$ is 2 too.\\
And so $f(x)=2x$, which indeed satisfies the given relation.

\pagebreak
\phantomsection\label{test2}{\color{blue}\section{Supremum and Infimum}}
\begin{tcolorbox}[colback=green!5!white,colframe=green!75!black]
Warm-Up Problem 3.1
  \tcblower
Does there exist a function $f:\mathbb{N}\to\mathbb{N}$ such that: $f(n)=f(f(n-1))+f(f(n+1))$ for all $n\in\mathbb{N}$?
\end{tcolorbox}
\textit{Solution:}
On the $RHS$ we have two terms of $f$ adding up to the one term on the $LHS$. This leads us to suspect that the range of $f$ should "flood" in $\mathbb{N}$. We formalize this insight by considering the $\displaystyle k=\min_{x>1, x\in\mathbb{N}} f(x)$, and let $n=\underset{x}{\operatorname{argmax}}f(x)\geq 2$.\\

Note that $k\geq 2$, otherwise $1=f(n)=f(f(n-1))+f(f(n+1))\geq 2$, absurd. Now, $f(n+1)\geq k\geq 2$, hence $k\leq f(f(n+1))=f(n)-f(f(n-1))<k$, a contradiction. Hence no such function exists.\\

When the range of a function is discrete and lower/upper bounded, of course, we can consider its minimum/maximum element, respectively. Now, let us define the notion of supremum and infimum, which generalizes this strategy/construction to the real numbers.\\

\begin{tcolorbox}[colback=red!5!white,colframe=red!75!black]
\textbf{Definition 2.2}: The infimum of a subset $S$ of a set $T$ is the greatest element in $T$ that is less than or equal to all elements of $S$, if such an element exists.\\
\textbf{Definition 2.3}: The supremum of a subset $S$ of a set $T$ is the least element in $T$ that is greater than or equal to all elements of $S$, if such an element exists.
\end{tcolorbox}
Note that an upper bounded set of real numbers has a supremum, and a lower bounded set of real numbers has an infimum.\\ 
\indent We also prove a lemma that we will use extensively later in this section:\\
\\
\begin{tcolorbox}[colback=blue!5!white,colframe=blue!75!black,title="Lemma"]
Assume that a function $f$ defined on an (non-trivial) interval $A$ of $\mathbb{R}$ has supremum $S$. Then there is a sequence $(a_{N})\subset{A}$ such that $\lim_{n\to\infty}f(a_{n})=S$.
\end{tcolorbox}
\textit{Proof:} Define the sequence $\epsilon_{1}=1, \epsilon_{n}=\frac{\epsilon_{n-1}}{10}$. There is a real number $a_{n}\in{A}$ such that $f(a_{n})>S-\epsilon_{n}$(otherwise, $S-\epsilon_{n}<S$ and $f(x)\leq S-\epsilon_{n}$ for all $x\in A$ contrary to the definition of the supremum) or $\epsilon_{n}>S-f(a_{n})\geq 0$ ($f(x)\leq S$ for all $x\in A$). \\
\indent As $n$ increases, $\epsilon_{n}\to 0$, hence, for all $\epsilon$ there is a $N$ such that $\epsilon_{n}<\epsilon$ for all $n\geq N$. By the construction of $(a_{n})$ it follows that $0\leq S-f(a_{n})\leq \epsilon_{n}\leq\epsilon$ for all $n\geq N$, hence, by the definition of the limit, $\displaystyle\lim_{n\to\infty}f(a_{n})=S$.  
\\
Similarly, we can also prove that the property holds for the infimum.
\\
\\
\indent When considering the supremum and infimum, it is often wise to prove $\inf f=\sup f$ in some range, to conclude that it is constant.\\
\begin{tcolorbox}[colback=green!5!white,colframe=green!75!black]
Problem 3.2
  \tcblower
   Let ${\varnothing \neq A \subset \mathbb{R}}$. Find all bounded functions ${g:A \to [1,+\infty)}$ such that $\displaystyle{\sup\left(\frac {f}{g}\right)=\frac {\sup f}{\inf g}}$ for all functions ${f:A \to [1,+\infty)}$ with $f\neq g$.
\end{tcolorbox}
\textit{Proof:}
Let $\displaystyle{ m=\inf \ g , \  M=\sup \ g. }$ If g is not constant, picking ${f=g^2}$ we obtain that $\displaystyle{mM=\sup \ g^2=M^2}$ (which is self-explained).
Thus, $\displaystyle{m=M}$ a contradiction, since we assume that $g$ is not constant. Hence g is constant.\\
\\
\begin{tcolorbox}[colback=green!5!white,colframe=green!75!black]
Problem 3.3
  \tcblower
   Find all functions $f:\mathbb{R^{+}}\rightarrow \mathbb{R^{+}}$ such that: $f(xf(y)+y)=f(xy)+f(y)$ for all $x,y\in\mathbb{R^{+}}$.
\end{tcolorbox}
\textit{Solution:} Denote by P(x,y) the given relation. If there is $a$ such that $f(a)<a$, then by $P(\frac{a}{a-f(a),a})$ yields $f(a)=0$, a contradiction. Hence $f(x)\geq x$ for all $x$.\\
\indent Now, $\displaystyle P(\frac{x}{y},y)\implies\displaystyle f(\frac{xf(y)}{y}+y)=f(x)+f(y)\geq \frac{xf(y)}{y}+y\hfill (1)$\\by the previous property.\\
\indent The relation (1) rewrites as: \[\frac{f(y)}{x}-\frac{y}{x}\geq \frac{f(y)}{y}-\frac{f(x)}{x}\]\hfill(2)\\
\indent The pith of the problem lies in proving that $f(1)=1$. Assume the opposite, i.e. that $C=f(1)\neq 1$. Since $f(1)\geq 1$, we deduce that $C=f(1)>1$, or $C-1>0$.\\
$P(1,y)\implies$ $f(f(y)+y)=2f(y)$.\\
\\
\begin{tcolorbox}[colback=olive!5!white,colframe=olive!75!black]
\indent \textbf{Claim:} There is some M such that $\displaystyle\frac{f(x)}{x}>C$ for all $x>M$
\end{tcolorbox}
\indent \begin{proof}Assume that there is no such $M$ such that $\displaystyle\frac{f(x)}{x}>C$ holds for all $x>M$. Then, we can construct an increasing, non-bounded, sequence of infinitely many terms $(x_{n})$ of positive reals such that $\displaystyle\frac{f(x_{i})}{x_{i}}\leq C$. By $\textbf{Bolzano Weirstrass Theorem}$ we may consider a sub-sequence of $x_{i}$, denote it by $a_{n}$ such that $\displaystyle\frac{f(a_{n})}{a_{n}}$ is convergent.\\
\indent Fix $y$ in (2) and set $x=a_{n}$ where $n$ runs through $\mathbb{N}$. This yields:\\
$\displaystyle\lim_{n\to\infty}(\frac{f(y)}{a_{n}}-\frac{y}{a_{n}}+\frac{f(a_{n})}{a_{n}})=\lim_{n\to\infty}\frac{f(a_{n})}{a_{n}}>\frac{f(y)}{y}$ hence \[C>\frac{f(y)}{y}\quad\text{for all}\;y\]\\
\indent Thus, $\displaystyle 1 \leq \frac{f(xf(y)+y)}{xf(y)+y}=\frac{f(xy)+f(y)}{xf(y)+y}<C$, where the equality followed from $P(x,y)$. Plugging in $y=1$ in the last one, we obtain:\\
$\displaystyle 1\leq \frac{f(x)+C}{xC+1}<C$, or $\displaystyle C+\frac{1}{x}-\frac{C}{x}<\frac{f(x)}{x}<C$. Tending $x\to \infty$, by the squeezing theorem, \[\displaystyle\lim_{x\to \infty}\frac{f(x)}{x}=C\]

Since $f(y)+y$ grows without bound as $y\to\infty$, we have \[C=\lim_{y\to\infty}\frac{f(f(y)+y)}{f(y)+y}=\lim_{y\to\infty}\frac{2f(y)}{f(y)+y}=\lim_{y\to\infty}\frac{\frac{2f(y)}{y}}{\frac{f(y)}{y}+1}=\frac{2C}{C+1}\] hence $C=1$ or $C=0$, a contradiction in both cases (since $C\neq 1$).\\
\end{proof}
Define $\displaystyle D=\inf_{x>M}\frac{f(x)}{x}$. Notice that $D\geq C>1$. Now, $P(1,y)$, with $y>M$ gives $f(f(y)+y)=2f(y)>D(f(y)+y)$ 
or $f(y)(2-D)>Dy$. Since $Dy>0$ and $f(y)>0$, we infer $2-D>0$ hence \[\frac{f(y)}{y}>\frac{D}{2-D}\quad\text{for all}\; y>M\]
Notice that $\displaystyle D<\frac{D}{2-D}$ (it boids down to $D>1$ which holds), which contradicts the definition of the infimum\hyperref[test01]{\color{blue}$^{2}$}.\\
\indent We have now established that $C=1$. To finish, we just bash:\\
$P(x.1)\implies$ $f(x+1)=f(x)+1$\hfill (3)\\
$P(\frac{1}{x},x)\implies $ $\displaystyle f(\frac{f(x)}{x}+x)=1+f(x)\geq \frac{f(x)}{x}+x$ hence \[(x-1)(f(x)-x)\geq 0\forall x\in\mathbb{R^{+}}\]\hfill (4)\\
\indent Picking in (4) $x<1$ we infer that $f(x)\leq x$ and since $f(x)\geq x$ for all $x$, we infer that $f(x)=x$ for all $x\in (0,1)$ and (2) implies that $f(x)=x$ for all $x$, as desired.\\
\begin{tcolorbox}[enhanced,fit to height=10cm,
  colback=green!25!black!10!white,colframe=green!55!black,title=\phantomsection\label{test01}{\color{white}Remarks},
  drop fuzzy shadow,watermark color=white]
\textit{Remark 1:} In the previous proof we used the definition of the supremum or infimum, but, essentially, we established lower bounds that we made them stronger by 'iterating' the given functional equation. We could have proven that $C=1$ by more elementary means presented in the previous section as follows:\\
Now, $P(1,y)$, with $y>M$ gives $f(f(y)+y)=2f(y)$. Applying the previous claim, $2f(y)=f(f(y)+y)>Cf(y)+Cy$ or $(2-C)f(y)>Cy$ for all $y>M$.\\
Since $Cy>0$ and $f(y)>0$, we infer $2-C>0$ hence $f(y)>\frac{C}{2-C}$. We need to formalize this pattern as explained in the previous section.\\
Consider the sequence $\displaystyle a_{n+1}=\frac{a_{n}}{2-a_{n}}$, where $a_{1}$. By induction on $n$ we infer that $f(y)>a_{n}y$ for all $n$.\\
Furthermore, we claim that $(a_{n})$ is strictly increasing. Indeed, $a_{1}=C>1$ and $\displaystyle a_{2}=\frac{C}{2-C}$ and we can check that $a_{2}>a_{1}$ is equivalent to $a_{1}=C>1$. Assume that $a_{n}>a_{n-1}$ for all $n<k$, hence $a_{k}>...>a_{1}=C>1$. We can check that $\displaystyle a_{k}<\frac{a_{k}}{2-a_{k}}=a_{k+1}$ is equivalent to $a_{k}>1$ as desired.\\
Since $(a_{n})$ is strictly increasing and upper-bounded ($a_{n}<2$), it converges, and we can calculate the limit of $(a_{n})$ by the last property of converging sequences. It turns out to be 1.\\
Hence, $a_{n}>...>a_{1}=C$ hence $1=\displaystyle\lim_{n\to\infty}a_{n}\geq C$, a contradiction.\\
\textit{Remark 2:} As the reader might have already observed, we can harness the existence of the supremum and/or the infimum on $\displaystyle\frac{f(x)}{x}$, to surpass the definition of the sequences in the 'Squeezing technique' explained in the previous section.
\end{tcolorbox}
\begin{tcolorbox}[colback=green!5!white,colframe=green!75!black]
Problem 3.4
  \tcblower
(AOPS) Determine all functions $f :\mathbb{R^{+}} \to \mathbb{R^{+}}$ such that $f(xy+f(x))=\frac{1}{2}\left(f(x)+f(y)\right)$, for all $x,y \in \mathbb{R^{+}}.$
\end{tcolorbox}
\textit{Solution:} Denote by $P(x,y)$ the given relation.\\
\indent First of all, we prove that $f$ is upper bounded. Set $C=f(1)$.\\
$P(1,y)$ gives $\displaystyle f(y+C)-C=\frac{f(y)-C}{2}$ (1). Consider $a_{1}=C>u>0$, $a_{n}=a_{n-1}+C$. Obviously, $(a_{n})$ diverges and by (1), $f(a_{n})-C=\frac{f(a_{n-1})-C}{2}$, hence $\lim_{n\to \infty}{(f(a_{n})-C)}=0$ or $\displaystyle\lim_{n\to \infty}f(a_{n})=C$.\\

Fix x, and for sufficiently large n, consider $\displaystyle P(x, \frac{a_{n}-f(x)}{x})$ to infer that $\displaystyle f(a_{n})=\frac{f(x)+f(\frac{a_{n}-f(x)}{x})}{2}>\frac{f(x)}{2}$. Tending $n\to\infty$, we obtain $f(x)\leq 2C$, i.e. f is upper-bounded.\\

Since $0<f(x)<2C$, consider $\displaystyle N=\inf_{x>0} f(x)$, $M=\displaystyle\sup_{x>0} f(x)$. We now bootstrap the sup and the inf:\\

Select two sequences $a_{n}$ and $b_{n}$ such that $\displaystyle\lim_{n\to\infty}f(a_{n})=M$ and \\$\displaystyle\lim_{n\to\infty}f(b_{n})=N$. 

Now, $\displaystyle P(b_{n}, \frac{y-f(b_{n})}{b_{n}})$ for all $\displaystyle y>N+1$, yields:\\
\[ f(y)=\lim_{n\to\infty}\frac{f(\frac{y-f(b_{n})}{b_{n}})+f(b_{n})}{2}\leq \frac{M+N}{2},\] hence $\displaystyle\sup_{y>N+1}f(y)\leq \frac{M+N}{2}$ (1).

Similarly, $P(a_{n}, \frac{y-f(a_{n})}{a_{n}})$ for all $y>M+1$, yields:\\
\[f(y)=\lim_{n\to\infty}\frac{f(\frac{y-f(a_{n})}{a_{n}})+f(a_{n})}{2}\geq \frac{N+M}{2},\] $\displaystyle\inf_{y>M+1}f(y)\geq \frac{M+N}{2}$ (2).\\

By $(1),(2)$ we deduce that $\displaystyle\sup_{y>\max(M+1,N+1)}f(x)=\inf_{y>\max(M+1,N+1)}f(x)$, hence $f$ is constant in $(\max(M+1,N+1), \infty)$, say $f(x)=D$.\\

Fix $x$ and select $y$ sufficiently large. $P(x,y)\implies$ $f(x)=D$, hence the only functions that satisfy $P(x,y)$ are the constant ones.\\
\\

We continue this section with an "almost-differential" equation that has a marvelous elementary solution that is in line with the ideas of this section.\\
\begin{tcolorbox}[colback=green!5!white,colframe=green!75!black]
Problem 3.5
  \tcblower
Determine all functions $f :\mathbb{R} \to \mathbb{R}$ such that $f(x)=f(\frac{x}{2})+f'(x)\frac{x}{2}$ $\forall x\in\mathbb{R}$.
\end{tcolorbox}

\textit{Solution: }For $x\neq 0 $ the given relation can be rewritten as: \[f'(x)=\frac{f(x)-f(\frac{x}{2})}{\frac{x}{2}}\] \hfill (1) \\Let $x_{0}>0$ be fixed. For that particular $x_{0}$ define the set: \[\text{A}=\{a\in (0,x_{0}] : f'(a)=f'(x_{0})\}\]

Since $\text{A}\subseteq [0,x_{0})$, then $\text{A}$ is bounded, hence $I=\inf \text{A}\geq 0$ is well-defined, and there exists a sequence $(a_{n})_{n=1}^{\infty}\in \text{A}$ such that $\displaystyle\lim_{n\to\infty}a_{n}=I$.\\

Since $f'(x)$ is continuous (owing to the continuity of the RHS of (1)), we obtain that $\displaystyle f'(x_{0})=\lim_{n\to\infty}f'(a_{n})=f'(\lim_{n\to\infty}a_{n})=f'(I)$. Hence, $I\in\text{A}$.\\

Assume, now, for contradiction that $I\neq 0$. Then by (1) and the mean-value theorem we conclude: $\displaystyle f'(x_{0})=f'(I)=\frac{f(I)-f(\frac{I}{2})}{\frac{I}{2}}=f'(\theta)$, where $\theta\in(\frac{I}{2},I)$. Hence, $\theta \in \text{A}$ and $\theta<I$, which is contrary to the definition of the infimum.\\

Therefore $I=0$, and as a result $f'(x_{0})=f'(I)=f'(0)$. Following a similar logic for negative values of $x_{0}$ (with the slight modification of considering the $\sup$ instead of the $\inf$ of the respective set) we may conclude that for every $x_{0}\in\mathbb{R}$ the relation $f'(x_{0})=f'(0)$ holds, i.e. $f'(x)=f'(0)=c$ $\forall x\in \mathbb{R}$, thus $f(x)\equiv cx+d\forall x\in \mathbb{R}$, as desired. We can check that all linear functions satisfy the given functional equation.\\
\begin{tcolorbox}[colback=green!5!white,colframe=green!75!black]
Problem 3.6
  \tcblower
Determine all functions $f :\mathbb{R^{+}} \to \mathbb{R}$ such that  $f(x + y) + f(x + z) - f(x)f(y + z) \geq  1$.
\end{tcolorbox}
\textit{Solution: }Denote by $P(x,y,z)$ the given relation. $P(2t,t,t)\implies$ $2f(3t)\geq 1+f(2t)^{2}\geq 1$, hence $f(x)\geq \frac{1}{2}$, thus $f$ is lower-bounded.\\
\indent Now, $P(t,t,t)\implies$ $(2-f(t))f(2t)\geq 1>0$ (1). Since $f(x)\geq\frac{1}{2}>0$, by (1) we infer that $2-f(x)>0\forall x$, hence $f$ is upper-bounded.\\
\indent Denote $S=\sup f(\mathbb{R^{+}})$, $I=\inf f(\mathbb{R^{+}})\geq \frac{1}{2}>0$. Now, $P(x,y,y)$ yields:\\
$2f(x+y)\geq f(x)f(2y)+1\geq I^{2}+1$ hence $f(x)\geq \frac{I^{2}+1}{2}$. However, by the definition of $I$, we should have: $I\geq \frac{I^{2}+1}{2}$, i.e. $I=1$.\\
\indent Pick a sequence $(a_{n})$ such that $\displaystyle\lim_{n\to\infty}f(a_{n})=S$. Consider $P(a_{n},y,z)$, which yields $2S\geq f(a_{n}+y)+f(a_{n}+z)\geq f(a_{n})f(y+z)+1$, therefore $2S\geq f(a_{n})f(y+z)+1$, hence $2S\geq \displaystyle\lim_{n\to\infty}f(a_{n})f(y+z)+1=Sf(y+z)+1$, hence $f(x)\leq \frac{2S-1}{S}$. However, by the definition of $\sup$, we have $S\leq \frac{2S-1}{S}$ hence $S=1.$\\
\indent All in all, $S=I=1$, i.e. $f(x)=1\forall x\in\mathbb{R^{+}}$, which works.\\ 
\begin{tcolorbox}[colback=green!5!white,colframe=green!75!black]
Problem 3.7
  \tcblower
Find all continuous functions $\displaystyle{f: \mathbb{R} \to (0,+\infty)}$ such that:
\[{\ln (f(f(x)))=f(x)}\quad\text{(1)},\; \displaystyle{f(-x)=\frac {1}{f(x)}} \quad\text{(2)} \forall \displaystyle{x \in \mathbb{R}}\]
\end{tcolorbox}
\textit{Solution: } In (2), x=0 gives that $f(0)=1=e^{0}$.\\
\indent Equation (1) gives $f(t)=e^{t}\forall t\in f(\mathbb{R})$. If we manage to prove that $f(\mathbb{R})=(0,\infty)$ then $f(t)=e^{t}\forall t>0$ and $f(-t)=\frac{1}{f(t)}=\frac{1}{e^{t}}=e^{-t}\forall t>0$, i.e. $f(x)=e^{x}\forall x\in \mathbb{R}$.\\
\indent Since f is continuous, $f(\mathbb{R})$ is an interval $\Delta$.

If $\Delta$ is upper-bounded, let $\sup\Delta =b$. Then $\displaystyle{f(b)=\lim_{x\to b^{-}}f(x)=\lim_{x\to b} e^x=e^b}$.  It is a fact that $f(b)\in \Delta$, that is, hence $\displaystyle{e^b\leq b,}$ a contradiction.

Let $\displaystyle{\inf\Delta =b\geq 0.}$ Then, since $\displaystyle{\lim_{x\to -\infty} f(x)= \lim_{x\to \infty}\frac{1}{f(x)}=\lim_{x\to \infty} e^{-x}=0 ,}$ we have $b\leq 0$, hence $\displaystyle{b=0}$.\\

All in all, $\displaystyle{\Delta=(0,+\infty)}$, qed.\\

We end this section with a very difficult problem that has recently appeared in an online competition. The following solution is the official solution which fits to the context of this section:\\
\begin{tcolorbox}[colback=green!5!white,colframe=green!75!black]
Problem 3.8
  \tcblower
\textbf (GIMO 2021 - Aops Mock contest) Determine all functions $f\colon\mathbb{R^{+}}\to\mathbb{R^{+}}$ such that:
$f(x)f(x + 2f(y)) = xf(x + y) + f(x)f(y)$, for all $x,y\in\mathbb{R^{+}}$.\\
\end{tcolorbox}
\textit{Solution:} Denote by $P(x,y)$ the given relation. Making arguments equal to cancel each other out, if there is an $a $ such that $f(a)<\frac{a}{2}$, $P(a-2f(a),a)$ yields $(a-2f(a))f(2a-2f(a))=0$ hence $f(a)=\frac{a}{2}$, a contradiction. If $f(a)=\frac{a}{2}$, $P(a,a)$ yields a contradiction, hence $f(x)>\frac{x}{2}$ for all $x\in\mathbb{R^{+}}$.\\
\indent Furthermore, note that (by the last established property) $f(x)f(x+2f(y))>\frac{xf(x)}{2}+f(x)f(y)$. Thus, $P(x,y)$ implies that $f(x+y)>\frac{f(x)}{2}$, hence whenever $x<y\implies f(x)<2f(y)$.\\
\begin{tcolorbox}[colback=olive!5!white,colframe=olive!75!black]
\indent\textbf{Claim 1:}$f(x)<2x+c\forall x\in\mathbb{R^+}$, for some constant $c$.
\end{tcolorbox}
\begin{proof}
Note that $x+y<x+2f(y)\implies f(x)f(x+2f(y))<2xf(x+2f(y))+f(x)f(y)$.

Fix an $a$, if $f(a)\le 2a$ then we are done. Else, \[2f(a)f(y)>2(f(a)-2a)f(a+2f(y))>a(f(a)-2a)+2f(y)(f(a)-2a)\] $\implies 4af(y)>a(f(a)-2a)\implies f(a)-2a<4f(y)$. The claim thus follows.
\end{proof}
\textbf{Identity:} Let $g(x) := 2f(x)$ and $h(x) := \frac{x}{f(x)}$, then for any $n\in\mathbb{N}$, we have
$$f(x+g^n(y))=h(x)^nf(x+y)+\sum_{i=1}^n h(x)^{n-i}f(g^{i-1}(y))$$
This can be done by induction, where the base case is $P(x,y)$ with both sides divided by $f(x)$.\\
\begin{tcolorbox}[colback=olive!5!white,colframe=olive!75!black]
\textbf{Claim 2:} $\lim_{x\to\infty} h(x)=1$.
\end{tcolorbox}

\begin{proof} First, note that since $h(x)<2$ for all $x>0$, there must exist constants $k_n$ for all $n$ such that $\sum_{i=1}^n h(x)^{n-i}f(g^{i-1}(1))<k_n$ for any $x$ and $n$.

Now, fix an $\epsilon<1$, then we can choose an $n\in\mathbb{N}$ such that $\epsilon^n<\frac{1}{10}$. If for some $x$ we have $h(x)\le\epsilon$ then we get$$\frac{x+g^n(1)}{2}<f(x+g^n(1))<h(x)^nf(x+1)+k_n<\frac{2x+2+c}{10}+k_n$$which can be true only for small enough $x$. Thus, there exists $N_{\epsilon}$ such that $h(x)>\epsilon$ whenever $x>N_{\epsilon}$.

On the other hand, if we fix an $\epsilon>1$, then we can choose an $n\in\mathbb{N}$ such that $\epsilon^n>10$. If for some $x$ we have $h(x)\ge \epsilon$ then we get$$2x+2g^n(1)+c>f(x+g^n(1))>h(x)^nf(x+1)>5x+5$$which, again, can only be true when $x$ is small enough. Thus, there exists $N_{\epsilon}$ such that $h(x)<\epsilon$ whenever $x>N_{\epsilon}$.

From the two above parts, we can conclude the claim.
\end{proof}
\begin{tcolorbox}[colback=olive!5!white,colframe=olive!75!black]
\textbf{Claim 3:} $f(x)\equiv x$
\end{tcolorbox}

\begin{proof} Fix an $x$ and let $y\to\infty$. $P(x,y)$ then gives

\begin{align*}
\frac{x}{f(x)}&=\lim_{y\to\infty}\frac{f(x+2f(y))-f(y)}{f(x+y)}\\
&= \lim_{y\to\infty}\frac{f(x+2f(y))}{f(x+y)}-\lim_{y\to\infty}\frac{f(y)}{f(x+y)}\\
&= \lim_{y\to\infty}\frac{h(x+y)}{h(x+2f(y))}\lim_{y\to\infty}\frac{x+\frac{2y}{h(y)}}{x+y}-\lim_{y\to\infty}\frac{h(x+y)}{h(y)}\lim_{y\to\infty}\frac{y}{x+y}\\
&= 1\times\lim_{y\to\infty}\left(\frac{2}{h(y)}-\frac{\left(\frac{2}{h(y)}-1\right)x}{x+y}\right)-1\times 1\\
&=1\times 2-1 =1
\end{align*}
and thus, we are done.\\
\end{proof}

\pagebreak
\phantomsection\label{test3}{\color{blue}\section{Limits Chasing}}
\begin{tcolorbox}[enhanced,fit to height=3.5cm,
  colback=green!25!black!10!white,colframe=green!75!black,title=Limits Chasing Technique,
  drop fuzzy shadow,watermark color=white]
If the function we are looking for has domain and co-domain the $\mathbb{R^+}$ and a
relation we have found for f has a variable x ”free” (i.e., it is not wrapped by any
”f”), we can set x → 0 (limit) which can yield something nice. If f is
also continuous (rare), we can do this, even if x is wrapped by some f. Moreover, we may as well replace $f(x)$ by its limit at some point $x_{0}$ by tending $x\to x_{0}$, assuming that we have found it (either with the squeeze theorem or with the theorem that we will discuss later on in this section).
\end{tcolorbox}
\begin{tcolorbox}[colback=green!5!white,colframe=green!75!black]
Problem 4.1
  \tcblower
  Find all functions $f: \mathbb{R^{+}}\to \mathbb{R^{+}}$ such that: $f(x^{2}+y^{2})+2f(x)f(y)=(x+y)^{2}$, for all $x,y \in \mathbb{R^{+}}$.
\end{tcolorbox}
\textit{Solution:} (by Rafail Tsiamis)\\ \indent By $(x+y)^{2}\leq 2(x^{2}+y^{2})$ we get $f(x^{2}+y^{2})\leq f(x^{2}+y^{2})+2f(x)f(y)=(x+y)^{2}\leq 2(x^{2}+y^{2})$, hence \[f(x)<2x\quad\forall{x}\in\mathbb{R^{+}}\]\hfill(1)\\ \indent Now (1) and the initial equation yield Q(x,y): $f(x^2+y^2) = (x+y)^2-2f(x)f(y) > (x+y)^2-8xy = \\=x^2 - y(6x-y)$ for all $x,y$. Considering $Q(\sqrt{x-\epsilon},\epsilon),$ for a small $\epsilon$ gives $f(x)\geq x-\epsilon A$ where A is a function of x. Tending $\epsilon \to 0$ we get $f(x)\geq x$ for all x.\\ \indent Using $P(x,y)$ we get $(x+y)^2=f(x^2+y^2)+2f(x)f(y) \geq x^2+y^2+2xy=(x+y)^{2}$ implying that $f(x)=x$ which finishes the proof.\\

Sometimes it is not easy to find the limits to employ the technique mentioned previously. The following theorem comes in handy in some cases:\\
\newtcbtheorem[auto counter,number within=section]{theo}%
  {Theorem}{fonttitle=\bfseries\upshape, fontupper=\slshape,
     arc=0mm, colback=blue!5!white,colframe=blue!75!black}{theorem}

\begin{theo}{Limits and bounds}
I If $f:\mathbb{R^{+}}\to \mathbb{R}$ is a monotonic function, then for each $x_{0}\geq 0$ the limits $\lim_{x \to x_{0}+}f(x)$, $\lim_{x \to x_{0}-}f(x)$ exist, and it also holds that:\\
a)If f is increasing:\\
$f(x_{0})\leq \lim_{x \to x_{0}+}f(x)=\inf_{x>x_{0}}f(x)$ and $\sup_{x<x_{0}}=\lim_{x\to x_{0}-}f(x)\leq f(x_{0})$\\
b)If f is decreasing:
$\sup_{x>x_{0}}f(x) \leq \lim_{x \to x_{0}+}\leq f(x_{0})$ and $f(x_{0}) \leq \lim_{x \to x_{0}-}f(x)=\inf_{x<x_{0}}f(x)$.
\end{theo}

This theorem allows us to consider the (left or right) limits and imposes bounds on them. Therefore, not only may we calculate the limits at specific points but also exclude some of their possible values that don't obey the bounds imposed. On a side note, this theorem is a generalization of the ninth property of the limits of sequences, though here, we deal with limits in real numbers. 
\\
\begin{tcolorbox}[colback=green!5!white,colframe=green!75!black]
Problem 4.2
\tcblower
Find all functions $f:\mathbb{R^{+}}\to \mathbb{R^{+}}$ such that: $f(xy+f(x))=f(x)f(y)+x$ for all $x,y \in \mathbb{R^{+}}$.\\
\end{tcolorbox}
\textit{Solution 1:} (by Athanasios Kontogeorgis)
Denote by $P(x,y)$ the given relation.
The proof unveils in 5 steps:\\
\textbf{Step 1:} $f(x)\geq x$ for all $x>0$.
\indent \begin{proof}
If there is an $a>0$ such that $a>f(a)$, $P(a,\frac{a-f(a)}{a})$ yields $f(a)=f(a)f(\frac{a-f(a)}{a})+a>a$, contrary to our assumption that $a>f(a)$.
\end{proof}
\textbf{Step 2:} $f(x)<1$ for every $x<1$\\
\indent \begin{proof}
Assume that there is an $a<1$ such that $f(a)\geq 1$. $P(a,\frac{f(a)}{1-a})$ yields $f(\frac{f(a)}{1-a})(1-f(a))=a\leq 0$, a contradiction.  
\end{proof}
\textbf{Step 3:} f is (weakly) increasing
\indent \begin{proof}
The proof utilizes an important technique:\\ 
\indent Assume for contradiction that $0<x_{2}<x_{1}$ with $f(x_{1})<f(x_{2})$. We consider the (positive) "ratio of change":\\
$\lambda=\frac{f(x_{2})-f(x_{1})}{x_{1}-x_{2}}$ or $x_{1}\lambda+f(x_{1})=x_{2}\lambda+f(x_{2})$. By "f-ing" the latter, we obtain:\\
$f(x_{1}\lambda+f(x_{1}))=f(x_{2}\lambda+f(x_{2}))$ or $f(x_{1})f(\lambda)+x_{1}=f(x_{2}f(\lambda)+x_{2}$  or $f(\lambda)=\frac{x_{1}-x_{2}}{f(x_{2})-f(x_{1})}=\frac{1}{\lambda}$.
Since $f(\lambda)\geq \lambda$ we have $\lambda\leq 1$. If $\lambda<1$, by Step 2, we obtain $f(\lambda)<1$ or $\frac{1}{\lambda}<1$ or $\lambda>1$, a contradiction.\\
\indent Hence $\lambda=1$ or $f(1)=1$.\\
Now, $P(1,y)\implies$ $f(y+1)=f(y)+1$ and by iterating $y\to y+1$ we obtain $f(y+n)=f(y)+n$(1) for all $n\in \mathbb{N}$.\\
$P(x,y+1)$ yields $Q(x,y): f(xy+x+f(x))=f(x)f(y)+f(x)+x=f(x)+f(xy+f(x))$. By considering $Q(x,\frac{y-f(x)}{x})$ for $y>f(x)$ we obtain $f(y+x)=f(x)+f(y)$ for all $y>f(x)$.\\
\indent Having obtained the property of "restricted additivity", it is common to extend this via the relation (1). This is how we do it:\\
\indent Fix an arbitrary pair $(x,y)$:\\
\indent Pick an integer $n$ such that $y+n>f(x)$. Hence, $f(y+x)=f(y+n+x)-n=f(y+n)+f(x)-n=f(y)+f(x)$, hence $f$ is additive.\\
\indent Since $f(x+y)=f(x)+f(y)>f(x)$, $f$ is strictly increasing, contrary to our assumption that $f$ is not increasing.
\end{proof}
\textbf{Step 4:} $\displaystyle\lim_{x\to 0}f(x)=0$ 
\begin{proof}
Since $f$ is increasing, the limit $L=\displaystyle\lim_{x\to 0+}f(x)=\lim_{x\to 0}f(x)$ exists and is finite.\\
\indent By $P(x,y)$ and Step 1 we deduce that: $f(x)f(y)+x\geq xy+f(x)$ hence $f^{2}(x)+x\geq x^{2}+f(x)$. Tending $x\to 0$ we obtain that $L^{2}\geq L$ hence $L\leq 0$ or $L\geq 1$.\\
\indent If $L\geq 1$, then $f(x)\geq L=\displaystyle\inf_{x>0}f(x)\geq 1$, a contradiction from Step 2.\\ 
\indent Hence $L\leq 0$. However, $f(x)>0$, thus $L\geq 0$ which implies $L=0$, as desired. 
\end{proof}
\textbf{Step 5:} $f(x)=x$ for all $x>0$
\begin{proof}
Tending $y\to 0$ in $f(x)f(y)+x\geq xy+f(x)$ we obtain that $x\geq f(x)$. This, in conjunction with Step 1 yields that $f(x)=x$ for all $x\in \mathbb{R^{+}}$.  
\end{proof}
\textit{Solution 2:} As previously, $f(x)\geq x \forall x$.\\

$P(x,y)$ and the previous relation gives: $\displaystyle x+f(x)f(y) \geq xy +f(x) \implies f(y)-1 \geq \frac{x}{f(x)}(y-1)$ (1). Since
$\displaystyle f(x) \geq x \implies 1 \geq \frac{x}{f(x)}$ (2)

Assume that there is $\displaystyle x_0$ such that $\displaystyle \frac{x_0}{f(x_0)}=m<1 \implies 1-m>0$. Plugging in (1) $x=x_0$, yields:\\
$\displaystyle f(y) \geq my-m+1>1-m \implies 0< \frac{y}{f(y)}< \frac{y}{1-m}$.
Now, tending $y\to 0^{+}$:
$\displaystyle \lim_{x \rightarrow 0^{+}} \frac{y}{f(y)}=0$.
Tending $x\to 0^{+}$ in (1):\\
$f(y) \geq 1 , \forall y>0 \implies f(\frac{1}{2}) \geq 1$. By $P(\frac{1}{2},2f(\frac{1}{2}))$ we have:
$f(2f(\frac{1}{2}))=f(2f(\frac{1}{2}))f(\frac{1}{2})+\frac{1}{2} \geq f(2f(\frac{1}{2}))+ \frac{1}{2} \implies 0 \geq 1/2$, a contradiction.\\

Hence, $\frac{x}{f(x)} \geq 1 \forall x$, which, in conjunction with (2), yields $f(x)=x$, which satisfies $P(x,y)$.\\
\begin{tcolorbox}[colback=green!5!white,colframe=green!75!black]
Problem 4.3
  \tcblower
(ISL 2020 A8) Determine all functions $f :\mathbb{R^{+}}\to \mathbb{R^{+}}$ such that: $f(x+f(xy))+y=f(x)f(y)+1$ for all x,y$\in\mathbb{R^{+}}$.
\end{tcolorbox}
\textit{Solution:} Denote by $P(x,y)$ the given relation.
If there are $a,b$ such that $f(a)=f(b)$, $P(1,a)-P(1,b)$ give $a-b=0$ hence $a=b$, a contradiction. Hence, $f$ is injective.\\  
\indent We claim that $x_{2}>x_{1}\implies f(x_{2})>f(x_{1})$, i.e. f is strictly increasing. Assume that $f(x_{2})\leq f(x_{1})$, and since f is injective, $f(x_{2})<f(x_{1})$.
$P(\frac{x}{y},y)\implies$ $Q(x,y): f(\frac{x}{y}+f(x))+y=f(\frac{x}{y})f(y)+1$.
Consider $\lambda=\frac{x_{2}-x_{1}}{f(x_{1})-f(x_{2})}>0$, or, equivalently, $f(x_{1})+\frac{x_{1}}{\lambda}=f(x_{2})+\frac{x_{2}}{\lambda}$.\\
$Q(x_{1}, \lambda)$ and $Q(x_{2}, \lambda)$ give that $f(\frac{x_{1}}{\lambda})=f(\frac{x_{2}}{\lambda})$, a contradiction, since $f$ is injective.\\
\indent Having proven that $f$ is strictly increasing, the limit $\displaystyle\lim_{x\to 0+}f(x)=L$ exists, as well as at every other point, besides 0.\\
\indent On $P(x,y)$, tend $x\to 0$ and notice that $x+f(xy)\to 0+L=L$, $f(x)\to L$, hence:\\
\[\displaystyle\lim_{x\to 0}f(x+f(xy))+y=\lim_{x\to L}f(x)+y=C+y\] \hfill(1)\\
and \[\displaystyle\lim_{x\to 0}f(x)f(y)+1=Lf(y)+1\]\hfill (2) \\
\indent By (1), (2) and $P(x,y)$ we obtain that $Lf(y)+1=C+y$ for every $y\in\mathbb{R^{+}}$ hence $f$ is linear. We can check that only $f(x)=x+1$ works.\\
\begin{tcolorbox}[colback=green!5!white,colframe=green!75!black]
Problem 4.4
  \tcblower
(Metrix 2020 - Aops Mock contest) Find all functions $f:\mathbb{R^{+}}\rightarrow \mathbb{R^{+}}$ such that $$f(x^2+2f(xy))=xf(x+y)+f(xy)$$for all $x,y\in \mathbb{R^{+}}$
\end{tcolorbox}
Here is the official solution:\\
\indent Denote by $P(x,y)$ the given relation.\\
\begin{tcolorbox}[colback=olive!5!white,colframe=olive!75!black]
\indent\textbf{Claim 1:} $f(x)>\frac{x}{2}$.
\end{tcolorbox}

\begin{proof}Assume the contrary at $x=a$.

Consider the equation $x^2+2f(a)=x+\frac{a}{x}$. As the right side is not smaller than the left side when $x=1$, but smaller when $x\to\infty$, this equation must have a root larger than or equal to $1$. Let this root be $r$.
We can see that $P(r,\frac{a}{r})$ gives contradiction.
\end{proof}
\begin{tcolorbox}[colback=olive!5!white,colframe=olive!75!black]
\indent \textbf{Claim 2:} $f$ is non-decreasing.
\end{tcolorbox}
\begin{proof} Assume the contrary that $f(b)>f(c)$ for some pair $b<c$. Let $d=2f(b)-2f(c)$.

Consider the equation $\sqrt{x}+\frac{b}{\sqrt{x}}= \sqrt{x+d}+\frac{c}{\sqrt{x+d}}$. After squaring both sides and clearing some calculation, we can see that this equation has a positive root $k$ and it satisfies the equation $(c-b)k(2k+c+b+2d)=d(b^2-dk-k^2)$. In particular, we have $dk+k^2<b^2$.

Let $l= \sqrt{k}+\frac{b}{\sqrt{k}}$. Comparing $P(\sqrt{k},\frac{b}{\sqrt{k}})$ and $P(\sqrt{k+d},\frac{c}{\sqrt{k+d}})$ gives $(\sqrt{k+d}-\sqrt{k})f(l)=f(b)-f(c)\implies 2f(l)= \sqrt{k+d}+\sqrt{k}$. From claim 1, we would have $\sqrt{k+d}> \frac{b}{\sqrt{k}}$ or $dk+k^2>b^2$, a contradiction.
\end{proof}
Claim 2 implies that $g(x)=\displaystyle\lim_{a\to x^+} f(a)$ exists for all $x\ge 0$.

First, taking $x\to 0^+$ in the original equation gives $g(2g(0))=g(0)$. If $g(0)>0$, then $g(2g(0))\ge f(2g(0))>g(0)$, a contradiction. Thus, $g(0)=0$. Taking $y\to 0^+$ in the original equation then gives $g(x^2)=xg(x)$.

Now, taking $x\to \infty$ in $P(\frac1x,x)$ gives $g(f(1))-f(1)=\displaystyle\lim_{x\to \infty} \frac{f(x+\frac1x)}{x}$. Let this value be $t$, then $\displaystyle\lim_{x\to \infty} \frac{f(x)}{x}= \lim_{x\to \infty} \frac{f(x+\frac1x)}{x+\frac1x}= \lim_{x\to \infty} \frac{f(x+\frac1x)}{x}\times \lim_{x\to \infty} \frac{x}{x+\frac1x} = t$.

Since $\frac{f(x)}{x}\le \frac{g(x)}{x}\le \frac{f(x+1)}{x}$ and $\displaystyle\lim_{x\to \infty} \frac{f(x+1)}{x}= \lim_{x\to \infty} \frac{f(x+1)}{x+1}\times \lim_{x\to \infty} \frac{x+1}{x}=t$, we also have $\displaystyle\lim_{x\to \infty} \frac{g(x)}{x}=t$ as well.

Now, from $g(x^2)=xg(x)$, we would have $\frac{g(x)}{x}=\frac{g(x^2)}{x^2}$, and thus $\displaystyle\frac{g(x)}{x}=\lim_{n\to\infty} \frac{g(x^{2^n})}{x^{2^n}}=t$ for all $x>1$. As $f(x)\le g(x)$ and $g(x-\epsilon)\le f(x)$ for all positive $\epsilon <x$, we must have $f(x)=tx$ for all $x>1$ as well.

Finally, $P(1,x)$ implies that $f(x)\equiv x$. Done.
\\
\begin{tcolorbox}[colback=green!5!white,colframe=green!75!black]
Problem 4.5
\tcblower
(Danube 2019 P2 Seniors) Find all increasing functions $f: \mathbb{R} \to \mathbb{R}$ such that $f(f(x^2) + y + f(y)) = x^2 + 2f(y)$, for all $x,y \in \mathbb{R}$.\\
\end{tcolorbox}
\textit{Solution} (by Minos Margaritis): Denote by $P(x,y)$ the given relation.
\begin{tcolorbox}[colback=olive!5!white,colframe=olive!75!black]
\indent\textbf{Claim 1:} f is injective
\end{tcolorbox}
\indent \begin{proof}
Assume that there are $a,b$ with $a<b$ such that $f(a)=f(b)$
We consider the following cases:\\
\indent \textbf{Case 1:}. $0\leq a$.\\
Considering $P(\sqrt{a},y)-P(\sqrt{b},y)$ yields what we want.\\
\indent \textbf{Case2:}. $a<0<b$\\
Since f is increasing, we have $f(x)=f(b) \forall x\in (0,b)$, hence, considering $a=x_{0}\in (0,b)$ we reduct this case to the previous case.\\
\indent \textbf{Case 3:} $a<b\leq 0$\\
f is obviously not upper bounded (by taking x sufficiently large in P(x,y))\\
Plugging in $x$ fixed but sufficiently large (such that $f(x^2)+y+f(y)>0 \forall y\in [a,b]$) and varying y in [a,b], where f is constant, we obtain that f is constant (as y varies) in an interval on the right of zero, which is absurd by the previous cases.
Hence $f$ is $1-1$.\\
\end{proof}
\begin{tcolorbox}[colback=olive!5!white,colframe=olive!75!black]
\textbf{Claim 2:} f is not lower bounded.
\end{tcolorbox}
\begin{proof}
Assume that m=$\inf(f(x))$ exists.
Since $f$ is strictly increasing: $\displaystyle\lim_{x\rightarrow-\infty}f(x)=m$, hence if we consider the sequence $x_{n}$ as $x_{n}=n \forall n\in \mathbb{N}$ and the sequence $y_{n}$ as $y_{n}=-f(x_{n}^2)-n,\forall n\in \mathbb{N}$,
we deduce, $\displaystyle\lim_{n\rightarrow +\infty}(f(x_{n}^2)+y_{n}+f(y_{n}))=-\infty$, thus $\displaystyle\lim_{n\rightarrow +\infty}f(f(x_{n}^2)+y_{n}+f(y_{n}))=m$.
Moreover, $\displaystyle\lim_{n\rightarrow +\infty}f(f(x_{n}^2)+y_{n}+f(y_{n}))=\lim_{n\rightarrow +\infty}x_{n}^2+2f(y_{n})=+\infty$ which is absurd.
Hence, $f$ is not lower bounded.
\end{proof}
Therefore, $\exists y_{0},f(y_{0})<0 \stackrel{P(x,y_{0})}{\Rightarrow} \exists x_{0},f(x_{0})=0$.
We have $P(0,x_{0}):f(f(0)+x_{0})=f(x_{0})$ and since f is 1-1 we have $f(0)=0$.
Now $P(x,0)\implies$ $f(f(x^2))=x^2,\forall x$.
Since $f$ is strictly increasing and involutive in $\mathbb{R^{+}}$, $f(x)=x, \forall x\geq 0$.
\indent Hence, $P(x,y)$ reduces to $f(x^2+y+f(y))=x^2+2f(y)$.
For each y there is an x such that $x^2+y+f(y)>0$, hence $y+f(y)=2f(y)$ or $f(y)=y,\forall y$ which satisfies our equation.

\pagebreak

\phantomsection\label{test4}{\color{blue}\section{Differentiability}}
Often in FE's, proving that a function is differentiable, even at some point, can be cumbersome. The following theorem often comes in handy:\\

\begin{theo}{Lebesgue Differentiation Theorem [1]}
A Every strictly monotone function in $\mathbb{R^{+}}$ is differentiable in $\mathbb{R^{+}}/A$ where A is a set of zero measure. In other words, a strictly monotone function is differentiable almost everywhere.\\
\end{theo}
For our purposes, we don't have to be pedantic to explain what a measure of a set is. Usually, we need just one point of differentiability and the previous theorem ensures the existence of it.
\\

\begin{tcolorbox}[colback=green!5!white,colframe=green!75!black]
Problem 5.1
  \tcblower
(ISL 2005 A2) Find all functions $f\colon \mathbb{R^{+}} \to \mathbb{R^{+}}$:  $f(x)f(y)=2f(x+yf(x))$ for all $x,y\in\mathbb{R^{+}}$.
\end{tcolorbox}
\textit{Solution:} Denote by $P(x,y)$ the given relation. We claim that $f$ is weakly increasing:\\
If there exists a pair $(x_{2},x_{1})$ such that $x_{2}>x_{1}$ and $f(x_{1})>f(x_{2})$, consider $\lambda=\frac{f(x_{1})-f(x_{2})}{x_{2}-x_{1}}>0$, or $\lambda f(x_{1})+x_{1}=\lambda f(x_{2})+x_{2}$(1).\\ (1), $P(x_{1},\lambda)$ and $P(x_{2},\lambda)$ give $f(\lambda)f(x_{2})=f(\lambda)f(x_{1})$ hence $f(x_{1})=f(x_{2})$, contrary to our assumption that $f(x_{1})>f(x_{2})$.\\
\indent Hence, $x>y\implies f(x)\geq f(y)$.\\
\\
\indent \textbf{Case 1:} f is not strictly increasing:\\
\indent If there is a pair $a>b$ such that $f(a)=f(b)$ then $f$ must be constant in the interval $[b,a]$. Then, $P(b, \frac{a-b}{f(b)})$ gives that $f(\frac{a-b}{f(b)})=2$ hence, there is an $A$ such that $f(A)=2$.\\
\indent Now, $P(A,y)$ gives that $f(y)=f(A+yf(A))$ hence $f$ is constant in $[y,A+yf(A)]$ or reversely, for all $y$. Spanning $y$ in $\mathbb{R^{+}}$ we conclude that $f$  is constant in an interval $[M,\infty)$.\\
\indent Fix $x$ and pick $y$ sufficiently large. $P(x,y)$ ($f(y)=f(x+yf(x))$ since f is constant for sufficiently large $y$) yields $f(x)=2$, hence $f(x)=2$ for all $x$ which satisfies the initial functional equation.\\
\\
\indent \textbf{Case 2:} f is strictly increasing:\\
\indent By the $\textbf{Lebesgue Monotone Differentiation Theorem}$, there is a point $x_{0}$ such that the limit $\displaystyle l=\lim_{h\to 0}\frac{f(x+h)-f(x)}{h}$ exists and is finite.\\
\indent By $P(x_{0},y):$ $\displaystyle\lim_{y\to 0}\frac{f(yf(x_{0})+x_{0})-f(x_{0})}{yf(x_{0})}=l=\lim_{y\to 0}\frac{\frac{f(x_{0})f(y)}{2}-f(x_{0})}{yf(x_{0})}=\lim_{y\to 0}\frac{f(y)-2}{2y}$, hence $\displaystyle\lim_{y\to 0}\frac{f(y)-2}{2y}=l$ exists.\\
\\ 
\indent Now, by the previous relation, we can see that the right-limit $\lim_{h\to 0+}\frac{f(x_{0}+h)-f(x_{0})}{h}$ is independent of $x_{0}$, hence $f$ is right differentiable at every point and the right-derivative equals $l$ at every point.\\
\indent As for the left limit, consider $P(x_{0}-y, \frac{y}{f(x_{0}-y)})$ to deduce that $f(x_{0}-y)f(\frac{y}{f(x_{0}-y)})=2f(x_{0})$ (1). Therefore:\\
$\displaystyle\lim_{y\to 0}\frac{f(x_{0})-f(x_{0}-y)}{y}\stackrel{(1)}{=}\lim_{y\to 0+}\frac{f(\frac{y}{f(x_{0}-y)}-2)}{\frac{2y}{f(x_{0}-y)}}$. Note that:\\
$f(x_{0}-y)>f(\frac{x_{0}}{2})$ if $y<\frac{x_{0}}{2}$ hence $\displaystyle 0\leq \lim_{y\to 0}\frac{y}{f(x_{0}-y)}\leq \lim_{y\to 0}\frac{y}{f(\frac{x_{0}}{2})}=0$, thus $\displaystyle\lim_{y\to 0}\frac{f(x_{0})-f(x_{0}-y)}{y}=\lim_{y\to 0}\frac{f(y)-2}{2y}=l$.\\
\indent All in all, $f'(x)=l$ for all $x$ hence $f$ is linear and we can check that the only possible solution is $f(x)=2$.\\

\begin{tcolorbox}[colback=green!5!white,colframe=green!75!black]
Problem 5.2
  \tcblower
(BMO 2022 P3) Determine all functions $f:\mathbb{R^{+}}\to \mathbb{R^{+}}$ such that:  $f(y(f(x))^{3} + x) = x^{3}f(y) + f(x)$ for all $x,y\in\mathbb{R^{+}}$.
\end{tcolorbox}
\textit{Solution:} As usual, denote by $P(x,y)$ the given relation. Considering $P(x,\frac{y}{f(x)^{3}})\implies$ $f(y+x)=x^{3}f(\frac{y}{f(x)^{3}})+f(x)$ hence $f(x+y)>f(x)$, thus f is strictly increasing.\\
\indent Hence, by \textbf{Lebesgue Differentiation theorem}, there is a point $x_{0}\in \mathbb{R^{+}}$ such that the limit $\displaystyle l=\lim_{x\to x_{0}}\frac{f(x)-f(x_{0})}{x-x_{0}}$ exists and is finite. Now, $P(y,x_{0})$ where $y\to 0$ yields that:
$\displaystyle\lim_{y\to 0}\frac{f(yf(x_{0}^{3})+x_{0})-f(x_{0})}{yf(x_{0})^{3}}=l=\frac{x_{0}^{3}}{f(x_{0})^{3}}\lim_{y\to 0}\frac{f(y)}{y}$, hence \[\displaystyle\lim_{y\to 0}\frac{f(y)}{y}=L\] exists and is finite.\\
\indent We prove that f is continuous:\\
\indent As proven, $\displaystyle\lim_{y\to 0}\frac{f(y)}{y}=L$ is finite, hence $\displaystyle\lim_{y\to 0}f(y)=0$.
Thus, \[\displaystyle\lim_{y\to 0}f(yf(x_{0}^{3})+x_{0})-f(x_{0})=\lim_{y\to x_{0}+}f(y)-f(x_{0})=x_{0}^{3}\lim_{y\to 0}f(y)=0\] hence $\displaystyle\lim_{x\to x_{0}+}=f(x_{0})$, hence f is right-continuous.\\ \indent We prove that f is left continuous as well. To that end, fix $x_{0}\in \mathbb{R^{+}}$, then $P(x_{0}-y, \frac{y}{f(x_{0}-y)^{3}})$ with $y<\frac{x_{0}}{2}$ yields:\\
\[f(x_{0})-f(x_{0}-y)=(x_{0}-y)^{3}f(\frac{y}{f(x_{0}-y)^{3}})<(x_{0}-y)^{3}f(\frac{y}{f(\frac{x_{0}}{2})^{3}})\], since $f$ is strictly increasing.\\
Hence, tending $y\to 0$, since $\displaystyle\lim_{y\to0}f(\frac{y}{f(\frac{x_{0}}{2})^{3}})=\lim_{x\to 0}f(x)=0$ we obtain that $\displaystyle\lim_{x\to x_{0}-}f(x)=f(x_{0})$, hence $f$ is left-continuous as well, i.e. continuous.
Now there are two continuations we can finish:\\
\\
\\
\textit{Continuation 1:} We prove that f is differentiable as well.\\
\indent Let $u$ be an arbitrary point. Then: $\displaystyle\lim_{y\to 0}\frac{f(yf(u)^{3}+u)-f(u)}{yf(u)^{3}}=\lim_{y\to 0}\frac{u^{3}f(y)}{yf(u)^{3}}=\lim_{x\to 0}\frac{f(y)}{y}\frac{u^{3}}{f(u)}=L\frac{u^{3}}{f(u)}=A$, i.e. the limit exists, hence $f$ is right-differentiable.
We claim that $f$ is left-differentiable as well:\\
By considering $P(u-v, v/f(u-v)^{3})$ we obtain that:\\
\[\lim_{v\to 0}\frac{f(u)-f(u-v)}{v}=\lim_{v\to 0}\frac{(u-v)^{3}f(\frac{v}{f(u-v)})}{v}=\lim_{v\to 0}\frac{f(\frac{v}{f(u-v)})}{\frac{v}{f(u-v)}}\frac{(u-v)^{3}}{f(u-v)}.\] \\
\indent Since $f$ is continuous, \[\lim_{v\to 0}\frac{(u-v)^{3}}{f(u-v)}=\frac{u^{3}}{f(u)}\] and $\displaystyle\lim_{v\to 0}\frac{f(\frac{v}{f(u-v)})}{\frac{v}{f(u-v)}}=\lim_{y\to 0}\frac{f(y)}{y}=L$ ( $\displaystyle\lim_{v\to 0} \frac{v}{f(u-v)}=0$ since $\displaystyle\lim_{v\to 0}f(u-v)=f(u)$), hence $\displaystyle\lim_{v\to 0}\frac{f(u)-f(u-v)}{v}=\\=L\frac{u^{3}}{f(u)}=A$.\\
\indent Therefore, the right-hand and left-hand limit of the ratio of change of $f$ around $u$ exist and are equal, hence $f$ is differentiable.\\
\indent To finish, by $P(x,y)$ we have: 
$\displaystyle\frac{f(yf(x)^{3}+x)-f(x)}{yf(x)^{3}}=\frac{f(y)}{y}\frac{x^{3}}{f(x)^{3}}$. Tending $y\to 0$, we obtain that $\displaystyle f'(x)=L\frac{x^{3}}{f(x)^{3}}$ hence $f'(x)f(x)^{3}=Lx^{3}$ for all $x\in \mathbb{R^{+}}$. For concreteness, we manually solve the latest differential equation:\\
$\displaystyle f'(x)f(x)^{3}=[\frac{f(x)^{4}}{4}]'=[\frac{Lx^{4}}{4}]'$, hence there is a constant $c\in\mathbb{R}$ such that $\frac{f(x)^{4}}{4}=\frac{Lx^{4}}{4}+c$. Tending $x\to 0$, since $\displaystyle\lim_{x\to 0} f(x)=0$ and $\displaystyle\lim_{x\to 0}x^{4}=0$ yields $c=0$, hence $f(x)=x$ for all $x\in \mathbb{R^{+}}$, as desired.\\

\textit{Continuation 2:} A more natural ending would be by calculating $f(1)$:\\
\indent Observe that $x>1\implies f(x)\geq 1$. Indeed, if there is an $a>1$ such that $f(a)<1$, $P(a,\lambda=\frac{a}{1-f(a)^{3}})$ implies that $f(\lambda)=a^{3}f(\lambda)+f(a)>a^{3}f(\lambda)\implies$ $a<1$, a contradiction.
Since, $f(x)>1$ for $x>1$ and $\displaystyle\lim_{x\to 0}f(x)=0<1$ and f is continuous, by the intermediate value property, we obtain that there is a $z$ such that $f(z)=1$.\\
\indent $P(y,z)\implies$ $f(y+z)=1+z^{3}>f(y)$ (since f is strictly increasing), equivalently, \[z^{3}>1-\frac{1}{f(y)}\]\hfill (1)\\ \indent Since $f$ is not upper-bounded (by $P(x,y)$) and strictly increasing, we obtain that $\displaystyle\lim_{x\to\infty}f(x)=\infty$. Tending $y\to\infty$ in (1), we deduce that $\displaystyle z^{3}\geq \lim_{y\to\infty}1-\frac{1}{f(y)}=1$ hence $z\geq 1$ (2).
If $f(1)<1$, $P(1,\frac{1}{1-f(1)^{3}})$ yields $f(1)=0$, absurd. Hence, $f(1)\geq 1=f(z)$, from which we infer that $1\geq z$ (3).\\
(2),(3) yield that $z=1$, hence $f(1)=1$.\\
$P(1,y)$ yields that $f(y+1)=f(y)+1$ (4) hence, inductively, $f(y+n)=f(y)+n$ (5) for all $n\in\mathbb{N}$.\\
\indent Now, $P(x,y+1)$, along with (4), yield $f(f(x)^{3}+y+x)=f(x+y)+x^{3}$, hence $f(f(x)^{3}+y)=f(x+y)+x^{3}$ for all $y>x$. We further extend this property for all $x,y$:\\
\indent Fix an arbitrary pair $(x,y)$ and pick an integer $n$ such that $y+n>x$. Then:\\
$f(f(x)^{3}+y+n)=f(y+n)+x^{3}$ hence $f(f(x)^{3}+y)=f(y)+x^{3}$ (6) (by (5) applied on both hands).\\
\indent Plugging in $x\to \sqrt[3]{f(x)})$ in (6), we deduce that $f(f(\sqrt[3]{f(x)})^{3}+y)=f(y)+f(x)$(7). Swapping x,y, we deduce that $f(f(\sqrt[3]{f(x)})^{3}+y)= f(f(\sqrt[3]{f(y)})^{3}+x)$ and by injectivity, $f(f(\sqrt[3]{f(x)})^{3})=x+c$ (8) for all $x\in\mathbb{R^{+}}$. Tending in (8) $x\to 0$, since $\displaystyle\lim_{x\to 0}\sqrt[3]{f(x)}=0$ we have that $\displaystyle\lim_{x\to 0}f(f(\sqrt[3]{f(x)})^{3}=0$, hence $c=0$.\\ \indent Now, $(7),(8)$ yield $f(x+y)=f(x)+f(y)$ for all $x,y\in \mathbb{R^{+}}$, which is Cauchy in $x,y\in \mathbb{R^{+}}$. Since $f$ is lower bounded ($f(x)>0$) it must be linear, and we can check that only $f(x)=x$ satisfies the initial functional equation.\\
\begin{tcolorbox}[enhanced,fit to height=7cm,
  colback=green!25!black!10!white,colframe=green!55!black,title=Remarks,
  drop fuzzy shadow,watermark color=white]
\textit{Remark 1:} We can surpass the calculation of $c$ by using the Cauchy-Type Lemma (CTL) mentioned in pg.6 of [3]. Indeed, $f(x+y+c)=f(x)+f(y)\implies$ $f(x)=c(x)+d$ for all sufficiently large $x$ where $c(x)$ is a solution to Cauchy's functional Equation. Note that $f(x)>0$ hence $c(x)>-d$ hence $c$ is lower-bounded in some interval $[M,\infty)$ hence $c(x)$ is linear.\\   
\textit{Remark 2:} For the record, all Cauchy-Type functional equations in $\mathbb{R^{+}}$ have linear solutions. In fact, a stronger statement holds:\\
Let $f,g,h:\mathbb{R}^+\to\mathbb{R}^+$ be functions such that
$$f(g(x)+y)=h(x)+f(y)$$holds for all $x,y\in\mathbb{R}^+$. Then $\frac{g(x)}{h(x)}$ is a constant function.\\
For the proof, have a look \href{https://artofproblemsolving.com/community/q1h2399644p19686093}{here.}
\end{tcolorbox}
\phantomsection\label{test5}{\color{blue}\section{Problems for Practice}}
\textbf{Problem 1:} Find all functions $f: \mathbb{R^{+}} \to \mathbb{R^{+}}$ such that:\\
\[f(x+f(x+y))=f(2x)+y\] for all $x,y \in \mathbb{R^{+}}$.\\
\\
\indent \textbf{Problem 2:} Solve Problem 3.3 without using supremum/infimum.\\
\\
\indent \textbf{Problem 3:} Find all $f:\mathbb{Z^{+}} \rightarrow \mathbb{R^{+}}$ such that $f(m)\geq m$ and\\ \[f(m+n) \mid f(m)+f(n)\] for all $m,n\in \mathbb{Z^+}$.\\
\\
\indent \textbf{Problem 4:} Find all functions $f: \mathbb{R^{+}} \to \mathbb{R^{+}}$ such that
\[f(xf(y))+f(yf(z))+f(zf(x))=xy+yz+zx\] for all $x,y,z \in \mathbb{R^{+}}$.\\
\\
\indent  \textbf{Problem 5:} Find all functions $f \colon \mathbb{R} \to \mathbb{R}$ that satisfy the inequality
\[ f(y) - \left(\frac{z-y}{z-x} f(x) + \frac{y-x}{z-x}f(z)\right) \leq f\left(\frac{x+z}{2}\right) - \frac{f(x)+f(z)}{2}\] for all real numbers $x < y < z$.\\
\\
\pagebreak
\phantomsection\label{test6}{\color{blue}\section{References}}

1. Pang-Cheng W. (2018). \href{https://cdn.artofproblemsolving.com/attachments/e/1/4edd5cbbb44a9f08d959399a9270448f5501f5.pdf}{Functional Equation}.\\
\\
\indent 2. Chen E. (2016). \href{https://web.evanchen.cc/handouts/FuncEq-Intro/FuncEq-Intro.pdf}{Introduction to Functional Equations}\\
\\
\indent 3. Kontogeorgis A., Tsiamis R. (2019). \href{https://arxiv.org/abs/1901.11131}{The CDE Method, a Technique in Functional Equations}\\
\\
\indent 4. \href{https://pdfcoffee.com/169-functional-equations-pdf-free.html}{169 Functional Equations}\\
\\
\indent 5. Andreescu T., Boreico I., Mushkarov O., Nikolov N. (2012). \textit{Topics in Functional Equations, 2nd Edition}, XYZ Press, LLC\\
\\
\indent 6. Pavardi A. H. (2018). \textit{Functional Equations in Mathematical Olympiads (2017 - 2018): Problems and Solutions (Vol. I)}, Kindle Edition \\
\\
\indent 7. Venkatachala B.J. (2002). \textit{Functional Equations: A Problem Solving Approach}, PRISM BOOKS PVT LTD, Bangalore
\BgThispage
\end{document}